\documentclass[a4paper,10pt]{article}

\usepackage{titletoc}
\usepackage{fancyhdr}
\usepackage[utf8]{inputenc}
\usepackage{amssymb}
\usepackage{amsmath}
\usepackage{amsthm}
\usepackage{verbatim}
\usepackage{fullpage}
\usepackage[utf8]{inputenc}
\usepackage{setspace}
\usepackage{geometry}
\usepackage{setspace}
 \usepackage{amsmath,amscd}
\usepackage{a4wide}
\usepackage{color}
\numberwithin{equation}{section}

\usepackage{ifpdf}
\ifpdf
    \usepackage[pdftex]{graphicx}   
    \usepackage[pdftex,     
            plainpages=false,   
            breaklinks=true,    
            colorlinks=false,
            pdftitle=My Document
            pdfauthor=My Good Self
           ]{hyperref}
    \usepackage{thumbpdf}
\else
    \usepackage{graphicx}       
    \usepackage{hyperref} {colorlinks=false}      
\fi

\newcommand{\al}{\alpha}
\newcommand{\R}{\mathbb{R}}
\newcommand{\N}{\mathbb{N}}

\newcommand{\p}{\partial}

\def \qed {\hfill \vrule height6pt width 6pt depth 0pt}

\newtheorem{Remark}{Remark}[section]
\newtheorem{Theorem}{Theorem}[section]
\newtheorem{Lemma}{Lemma}[section]
\newtheorem{Definition}{Definition}[section]
\newtheorem{Assumption}{Assumption}[section]

\begin{document}

 \title {Mean Field Limit and Propagation of Chaos for a Pedestrian Flow Model}


\author {Li Chen\footnotemark[1]  \, \quad Simone Göttlich\footnotemark[1] \, \quad Qitao Yin\footnotemark[1]}

\footnotetext[1]{University of Mannheim, Department of Mathematics,
68131 Mannheim, Germany. E-mail: \{chen, goettlich\}@math.uni-mannheim.de; qyin@mail.uni-mannheim.de.}

 \maketitle

\begin{abstract}
In this paper a rigorous proof of the mean field
  limit for a pedestrian flow model in two dimensions is given by using a probabilistic method.
  The model under investigation is an interacting particle system coupled to the eikonal equation on the microscopic scale.
  For stochastic initial data, it is proved that the solution of the $N$-particle pedestrian flow system with properly chosen cut-off
  converges in the probability sense to the solution of the characteristics of the non-cut-off Vlasov equation.
  Furthermore, the result on propagation of chaos is also deduced in terms of bounded Lipschitz distance.
\end{abstract}
 {\bf Keywords:} probabilistic method, pedestrian flow, mean field limit, Vlasov equation, propagation of chaos.\\
{\bf AMS Classification:} 35Q83, 82C22
\section{Introduction}


The notable interest of pedestrian flow models can be dated back to four
decades ago with a considerable increase in interest since about year 2000.
For a general recent overview we refer
to~\cite{Bellomo2011,Bellomo2012,Colombo2008,Cristiani2011,Degond2013,HM95,Hughes2001,Piccoli2009}
and the references therein.
Pedestrian models share striking analogies in classical
physics such as gases and fluids, but are also applied to
the description of opinion formation \cite{Toscani2006}, group dynamics or other
social phenomena \cite{Naldi2010}.
Pedestrian flow models are an ideal starting point for the derivation of other or
more general quantitative behavioral models, since the relevant quantities
of pedestrian motions are easily measured so that corresponding models are comparable with
empirical data \cite{HM95}.
The modelling presented here is based on the idea that
behavioral changes are guided by so-called social fields or social forces, which
have been suggested by Lewin \cite{Lewin51}. Numerical simulations have been recently
carried out in \cite{EGKT14} on the microscopic and
macroscopic level using the finite particle method (FPM).
Some interesting spatiotemporal patterns are observed.

This paper provides the detailed derivation from the $N$-particle (pedestrian) Newtonian
system to its mean field limit or Vlasov equation.
Instead of the formal derivation with the help of the BBGKY hierarchy, which can be found in \cite{EGKT14, Sophn12}, we will
rigorously derive the kinetic description by a probabilistic method, which is inspired by Boers and Pickl \cite{BP15}, Hauray and Jabin
\cite{HJ07, HJ11}, Philipowski \cite{Philip2007} and Sznitman \cite{Sznitman91} and all the references therein.

However, the proposed pedestrian model involves a singularity, which
comes from the albeit bounded interaction force and is similar to the one generated by the Coulomb potential
in 2-$d$. While the authors in \cite{BP15} do not tackle the direct Coulomb potential in 3-$d$,
i.e. they consider the singularity that is a little weaker than for the Coulomb potential,
we are capable to deal with the singularity directly due to the compact support of the considered interaction force.
Another difficulty lies in the treatment of the dissipative terms since the interaction force
depends not only on the position $x$ but also on the velocity $v$. This will lead to
extra work on the estimates and is up to our knowledge rarely done before.

We now briefly explain our approach. In order to obtain the convergence
between the exact and the mean field dynamics, we mainly split the proof into two parts:
Using the Newtonian system with cut-off as a starting point,
we show that the Newtonian and the intermediate system (Vlasov flow with cut-off) are close to each other for $N$ being large enough.
The next step is to show the intermediate system converges to the Vlasov flow without cut-off.
Inbetween we use characteristics as a bridge to connect the Newtonian system and the mean field dynamics.

Additionally, assuming stochastic initial data offers a way to rule out those deterministic dynamics that do not
fit into the proper configuration of the Vlasov equation in the sense that those particles have small probability to appear.
In doing so, we obtain the convergence in (probability) measure between the exact and the mean field dynamics.
As a direct implication of the convergence, we prove the
propagation of chaos in terms of bounded Lipschitz distance.

This article is organized as follows: we start with the introduction
of the pedestrian flow model in Section~\ref{sec:model}.
Then, in Section~\ref{sec:prelim} some notations and preliminary work will be introduced.
In Section~\ref{sec:mean} we state the main results and present the corresponding proofs.
Section~\ref{sec:prog} is devoted to the propagation of chaos.
At this point, we also refer to \cite{Sznitman91} for other classical results with bounded Lipschitz continuity.
Finally, we summarize our results.


\section{Modeling of Pedestrian Flow} \label{sec:model}


Following the pedestrian flow model originally introduced in~\cite{EGKT14},
we consider a two-dimensional interacting particle system with position $x_i \in \R^2$
and velocity $v_i \in \R^2, i=1, \ldots, N$.
The equations of motion read
 \begin{align}
\label{NS} \begin{cases}
 \displaystyle{\frac{dx_{i}}{dt}}=v_{i}, \\ \vspace{0.005cm}\\
 \displaystyle{\frac{dv_{i}}{dt}}=\frac{1}{N-1}\sum_{i \neq j}F(x_i-x_j, v_i-v_j)+G(x_i,v_i),
 \end{cases}
 \end{align}
where $F(x,v)$ denotes the total interaction force and $G(x,v)$ the desired velocity
and direction acceleration. More precisely, $F(x,v)$ consists of the interaction force $F_{int}(x)$ and the dissipative force
$ F_{diss}(x,v)$, i.e.,
$$F(x,v)=(F_{int}(x)+F_{diss}(x,v))\mathcal{H}(x,v)$$
with
$$
  F_{int}(x) =k_n\frac{x}{|x|}(2R-|x|)
  = 2Rk_n \frac{x}{|x|}- k_nx,
$$
\begin{eqnarray*}
F_{diss}(x,v) &=& F^n_{diss}(x,v)+F^t_{diss}(x,v)\\
  &=&-\gamma_n \frac{\langle v,x \rangle}{|x|^2}x-\gamma_t\left(v-\frac{\langle v,x \rangle}{|x|^2}x\right)\\
  &=& \frac{\langle v,x \rangle}{|x|^2}(\gamma_t-\gamma_n)x-\gamma_tv
\end{eqnarray*}
and
$$\mathcal{H}(x,v):=\mathcal{H}_{2R}(|x|) \cdot \widetilde{\mathcal{H}}_{2\widetilde{R}}(|v|),$$
where $\mathcal{H}_{2R}(|x|)$ and $ \widetilde{\mathcal{H}}_{2\widetilde{R}}(|v|)$
are smooth functions with compact support that satisfy
$$
\mathcal{H}_{2R}(|x|)=\begin{cases}
0, & |x| >2R,\\
1, & |x| < R,
\end{cases}
\quad \hbox{and} \quad
 \widetilde{\mathcal{H}}_{2\widetilde{R}}(|v|)=\begin{cases}
 0, & |v| >2\widetilde{R},\\
 1, & |v| < \widetilde{R}.
\end{cases}
$$
Here, $F^n_{diss}(x,v)$ and $F^t_{diss}(x,v)$ are the normal
dissipative force and the tangential friction force, respectively.
Moreover, $k_n$ is the interaction constant and $\gamma_n, \gamma_t$ are
suitable positive friction constants.

\begin{Remark}
  To obtain a realistic behavior of pedestrians, the functions $\mathcal{H}_{2R}(|x|)$ and
	$\widetilde{\mathcal{H}}_{2\widetilde{R}}(|v|)$
	are used to express that the interaction force and the pedestrian velocity
  are of finite range. Mathematically, the total force is considered on a bounded domain.
\end{Remark}

\noindent The desired velocity and direction acceleration is given by
$$
G(x,v):=G(x,v, \rho)=\frac{1}{T}\left(-U(\rho)\frac{\nabla \Phi(x)}{|\nabla
\Phi(x)|}-v\right),
$$
where
$$\rho=\rho(x)=\frac{1}{N^R_{\max}} \sum_{j, |x-x_j|<R}1.$$
$N_{max}^R$ depends on the time $t$ via the coupling to the positions $x_j$. For a fixed time $t$, $N_{max}^R$ describes
the maximal number of particles in a ball of radius $R$ and is used here as a normalization parameter. This means,
we only scale the number of particles in this region in the sense how compressed they are.
$\Phi$ is given by the solution of the eikonal equation
$$
(\rho(x))|\nabla\Phi|-1=0,
$$
where $U: [0,1] \to [0, U_{\max}]$ is a density-dependent velocity function.
The reaction time $T$ might also depend on the density $\rho$.

The kinetic equation associated with this particle system describes the evolution
of the (effective one particle) density $f(t,x,v)$ as
\begin{eqnarray}
\label{VE}
\p_tf+v \cdot \nabla_xf+\nabla_v \cdot \left[(F*f)  f\right]+\nabla_v \cdot (G f)=0.
\end{eqnarray}
See~\cite{EGKT14} for more details and the derivation of macroscopic models
for different moment closures.


\section{Notations and Preliminary Work} \label{sec:prelim}


Now, we consider the pedestrian flow model~\eqref{NS} with cut-off of order $N^{-\theta}$ with $0< \theta < \frac{1}{4}$, i.e.,
$$F^N(x,v)=\begin{cases}
\left(2Rk_n \displaystyle{\frac{x}{|x|}}- k_nx+\displaystyle{\frac{\langle v,x \rangle}{|x|^2}}(\gamma_t-\gamma_n)x-\gamma_tv
\right) \mathcal{H}(x,v), & |x| \ge N^{-\theta}, \\
\left((2Rk_nN^{\theta}-k_n)x+N^{2\theta}\langle v,x \rangle (\gamma_t-\gamma_n)x-\gamma_tv\right) \mathcal{H}(x,v), & |x| < N^{-\theta}. \end{cases}$$

In order to present the analytical results in Section~\ref{sec:mean} in a concise and clear manner, we restrict to the following notations. 
\begin{Definition}
  \begin{enumerate}
    \item Let $(X^N_t, V^N_t)$ be the trajectory on $\R^{4N}$ which evolves according to the Newtonian
    equation of motion with cut-off, i.e.,
    \begin{align}
\label{NF} \begin{cases}
 \displaystyle \frac{d}{dt}X^N_t=V^N_{t}, \\ \vspace{0.005cm}\\
 \displaystyle{\frac{d}{dt}}V^N_t=\Psi^N(X^N_t, V^N_t)+\Gamma(X^N_t, V^N_t),
 \end{cases}
 \end{align}
    where $\Psi^N(X^N_t, V^N_t)$ denotes the total interaction force with
    $$\displaystyle \big(\Psi^N(X^N_t, V^N_t)\big)_i=\frac{1}{N-1}\sum_{i \neq j}F^N(x^N_i-x^N_j, v^N_i-v^N_j),$$
    while $\Gamma(X^N_t, V^N_t)$ stands for the desired velocity and direction
    acceleration with
    $$\big(\Gamma(X^N_t, V^N_t)\big)_i=G(x^N_i,v^N_i).$$

    \item Let $(\overline{X}^N_t, \overline{V}^N_t)$ be the trajectory on $\R^{4N}$ which evolves according to
    the Vlasov equation
    \begin{eqnarray} \label{vlasovN}
      \p_t f^N+v \cdot \nabla_x f^N+\nabla_v \cdot [(F^N*f^N) f^N]+\nabla_v \cdot (G
      f^N)=0,
    \end{eqnarray}

    i.e.,
    \begin{align}
\label{VF} \begin{cases}
 \displaystyle \frac{d}{dt}\overline{X}^N_t=\overline{V}^N_{t}, \\ \vspace{0.005cm}\\
 \displaystyle{\frac{d}{dt}}\overline{V}^N_t=\overline{\Psi}^N(\overline{X}^N_t, \overline{V}^N_t)+\Gamma(\overline{X}^N_t, \overline{V}^N_t),
 \end{cases}
 \end{align}
    where $\displaystyle \big(\overline{\Psi}^N(\overline{X}^N_t, \overline{V}^N_t)\big)_i=\iint F^N(\overline{x}^N_i-y, \overline{v}^N_i-w)f^N(t, y,w)\,dydw$ and
    $\big(\Gamma(\overline{X}^N_t,
    \overline{V}^N_t)\big)_i=G(\overline{x}^N_i,\overline{v}^N_i)$ represent the total
    interaction force and the desired velocity and direction
    acceleration, respectively.
    \end{enumerate}
\end{Definition}

\noindent  If $N$ is removed from the superscript, then $(X_t, V_t)$ and $(\overline{X}_t, \overline{V}_t)$
  denote the particle configurations driven by the force without cut-off.
  Analogically, if $t$ is removed from the subscript, $(X, V)$ and $(\overline{X}, \overline{V})$
  represent the stochastic initial data, which are independent and identically distributed.
 Note that we always consider the same initial data for both systems, that means $(X, V) = (\overline{X},
  \overline{V})$.

\begin{Remark}
We also point out several facts for the interaction force $F^N(x,v)$ with cut-off
and the acceleration $G(x,v)$. All the properties can be checked by direct computations.
  \begin{itemize}
    \item[(a)] $F^N(x,v)$ is bounded, i.e., $|F^N(x,v)| \le C$.
    \item[(b)] $F^N(x,v)$ satisfies the following property
    $$|F^N(x,v)-F^N(y,v)| \le q^N(x,v) |x-y|,$$
   where $q^N$ has compact support in $B_{2R} \times B_{2\widetilde{R}}$ with
            $$
            q^N(x,v):=
              \begin{cases} \displaystyle
                C \cdot \frac{1}{|x|}+C, & |x|\ge N^{-\theta}, \\ \vspace{0.005cm}\\
                C \cdot N^{\theta}, & |x|< N^{-\theta}.
              \end{cases}
              $$
    \item[(c)] $F^N(x,v)$ is Lipschitz continuous in $v$.
    \item[(d)] $G(x,v)$ is bounded, i.e., $|G(x,v)| \le C$.
  \end{itemize}
  In this context, we use $C$ as a universal constant that might depend on $k_n, R, \widetilde{R}, \gamma_n, \gamma_t$ .
\end{Remark}

	\noindent Furthermore, if there is a singularity in the velocity $v$ in the interaction
  potential similar to Remark 3.1(b), i.e.,
	
  $$|F^N(x,v)-F^N(x,w)| \le \widetilde{q}^N(x,v) |v-w|,$$
  where $\widetilde{q}^N(x,v)$ has compact support in $B_{2R} \times B_{2\widetilde{R}}$ with
  $$
            \widetilde{q}^N(x,v):=
              \begin{cases} \displaystyle
                C \cdot \frac{1}{|v|}+C, & |v|\ge N^{-\theta}, \\ \vspace{0.005cm}\\
                C \cdot N^{\theta}, & |v|< N^{-\theta},
              \end{cases}
              $$
it can be treated by using the same method as above and the results also
apply.
%


\section{Mean Field Limit} \label{sec:mean}


In this section, we present our key results in full detail. To show the desired convergence, our method can be summarized as follows.
First, we start from the Newtonian system with carefully chosen cut-off and meanwhile introduce an intermediate system which
involves convolution-type interaction with cut-off, namely \eqref{vlasovN} and \eqref{VF}. Then, we show the convergence of the intermediate system to
the final mean field limit, where the law of large number comes into play. The
crucial point of this method is that we apply stochastic initial data or in other
words we consider a stochastic process.
It enables us to use the tools from
probability theory, which helps to better understand the mean field process.
The overall procedure can be summarized as follows:

\begin{tiny}
\[
\begin{CD}
\begin{cases}
 \displaystyle \frac{dx_{i}}{dt}=v_{i} \\ \vspace{0.005cm}\\
 \displaystyle \frac{dv_{i}}{dt}=\frac{1}{N-1}\sum_{i \neq j}F(x_i-x_j, v_i-v_j)+G(x_i,v_i)
 \end{cases}
 @> \hbox{cut-off}>> \quad
\begin{cases}
 \displaystyle{\frac{dx^N_{i}}{dt}}=v^N_{i} \\ \vspace{0.005cm}\\
 \displaystyle{\frac{dv^N_{i}}{dt}}=\frac{1}{N-1}\sum_{i \neq j}F^N(x^N_i-x^N_j, v^N_i-v^N_j)+G(x^N_i,v^N_i)
 \end{cases}
\\
@. @VV N\to \infty V\\
\p_tf+v \cdot \nabla_xf+\nabla_v \cdot \left[(F*f)  f\right]+\nabla_v \cdot (G f)=0
 @ < \hbox{without cut-off}<< \quad
  \p_t f^N+v \cdot \nabla_x f^N+\nabla_v \cdot [(F^N*f^N) f^N]+\nabla_v \cdot (G
      f^N)=0\\
      @VV \hbox{characteristics} V @VV \hbox{characteristics} V\\
      \begin{cases}
 \displaystyle \frac{d \overline{x}_{i}}{dt}=\overline{v}_{i} \\ \vspace{0.005cm}\\
 \displaystyle \frac{d\overline{v}_{i}}{dt}=\iint F(\overline{x}_i-y, \overline{v}_i-w)f(t, y,w)\,dydw+G(\overline{x}_i,\overline{v}_i)
 \end{cases}
 @.
  \begin{cases}
 \displaystyle \frac{d\overline{x}^N_{i}}{dt}=\overline{v}^N_{i} \\ \vspace{0.005cm}\\
 \displaystyle \frac{d\overline{v}^N_{i}}{dt}=\iint F^N(\overline{x}^N_i-y, \overline{v}^N_i-w)f^N(t, y,w)\,dydw+G(\overline{x}^N_i,\overline{v}^N_i)
 \end{cases}
\end{CD}
\]
\end{tiny}

\noindent The following assumptions are used throughout this section.

\begin{Assumption}
  We assume that
  \begin{enumerate}

    \item[(a)] there exists a time $t>0$ and a constant $C$ such that the solution $f(t,x,v)$ of the Vlasov equation \eqref{VE}
    satisfies
    $$\sup_{0\le s \le t} \left|\left| \iint \frac{1}{|x-y|}f(s,y,v)\,dydv \right|\right|_{\infty} \le C,$$

    \item[(b)] the function $G(x,v)$ is Lipschitz continuous both in $x$ and $v$, i.e., there exists a constant $L$
    such that
$$|G(x,v)-G(x',v')| \le L \, (|x-x'|+|v-v'|), \quad \forall \, (x,v), (x',v') \in \R^{4N}.$$
  \end{enumerate}
\end{Assumption}

%

\begin{Definition}
  Let $\al \in (0, \frac{1}{5})$ and $S_t : \R^{4N} \times \R \to \R$ be the stochastic process given by
  $$S_t=\min\Big\{1, N^{\al}\sup_{0\le s \le t}\Big|(X^N_s,V^N_s)-(\overline{X}^N_s, \overline{V}^N_s)\Big|_{\infty}\Big\}.$$
  The set, where $|S_t|=1$, is defined as $\mathcal{N}_{\al}$, i.e.,
 \begin{eqnarray} \label{nalpha}
\mathcal{N}_{\al}:=\left\{(X,V) : \sup_{0\le s \le t}\Big|(X^N_s,V^N_s)-(\overline{X}^N_s, \overline{V}^N_s)\Big|_{\infty}> N^{-\al}\right\}.
 \end{eqnarray}
\end{Definition}

  \noindent Here and in the following we use $|\cdot|_{\infty}$ as the supremum norm on $\R^{4N}$.
  \noindent Note that
  $$\mathbb{E}_0(S_{t+dt}-S_{t} \,|\,\mathcal{N}_{\al}) \le 0,$$
  since $S_t$ takes the value of one for $(X,V) \in \mathcal{N}_{\al}$.

\begin{Theorem}
 Let $\theta \in (0, \frac{1}{4})$, $\al \in (0, \frac{1}{5})$, $\beta \in (\al, \frac{1-\alpha}{4})$, $\gamma \in (0, \frac{1-\al}{4}-\theta)$ and $f^N(t,x,v)$
 be the solution to the Vlasov equation \eqref{vlasovN}. Suppose that $f^N(t,x,v)$ satisfies Assumption 4.1(a) and Assumption 4.1(b) holds for $G(x,v)$.
 Then there exists a constant $C$  such that
  $$\mathbb{P}_0 \left(\sup_{0\le s \le t}\left|(X^N_s,V^N_s)-(\overline{X}^N_s, \overline{V}^N_s)\right|_{\infty}> N^{-\al} \right) \le e^{Ct}\cdot N^{-n},$$
  where $n=\min \{ 1-\al-4\beta, 1-\al-4\theta-4\gamma, \beta-\al \}$. Furthermore, if  $f^N(t,x,v) \in L^{\infty}((0,\infty); L^1(\R^2\times \R^2)) \cap
  L^{\infty}( (0,\infty);L^{\infty}(\R^2\times \R^2))$,  it holds with a $\theta$-independent convergence rate that
  $$\mathbb{P}_0 \left(\sup_{0\le s \le t}\left|(X^N_s,V^N_s)-(\overline{X}^N_s, \overline{V}^N_s)\right|_{\infty}> N^{-\al} \right)
  \le e^{Ct}\cdot r(N),$$
  where the convergence rate $r(N)=\max\{N^{
  -(1-\al-4\beta)},N^{\al-\beta},N^{-(1-\al-4\gamma)}\ln^2N\}$.
\end{Theorem}

 \noindent   We remark that $f^N(t,x,v) \in L^1 (\R^2\times \R^2)$ is automatically satisfied due to the mass conservation.\\

With additional assumption on the initial condition $f_0$ for the equations \eqref{VE}, \eqref{vlasovN} and on the solution of
the Vlasov equation without cut-off, we further extend our result to

\begin{Theorem}
  Let $f(t,x,v)$ and $f^N(t,x,v)$ be the solution to the Vlasov equation \eqref{VE} and \eqref{vlasovN} respectively with the same
  initial data $f_0$. Suppose that Assumption 4.1(b) is satisfied. Moreover, $\nabla f_0$ is integrable and $f(t,x,v) \in 
  L^{\infty}( (0,\infty);L^{\infty}(\R^2\times \R^2))$. Then there holds
  $$\lim_{N \to \infty}\mathbb{P}_0 \left(\sup_{0\le s \le t}\Big|(X^N_s,V^N_s)-(\overline{X}_s, \overline{V}_s)\Big|_{\infty}> N^{-\al} \right) =0.$$
\end{Theorem}

\noindent The proofs of both theorems will be presented at the end of this section.\\

  The additional requirement on $f(t,x,v)$ stems from the existence and uniqueness of the solution to the Vlasov
  equation, which will be shown in another independent work in the near future.

\begin{Definition} Let $\beta \in (\al, \frac{1-\al}{4}), \gamma \in (0, \frac{1-\al}{4}-\theta)$. The sets $\mathcal{N}_{\beta}$ and $ \mathcal{N}_{\gamma}$
  are characterized by
\begin{eqnarray} \label{nbeta}
  \mathcal{N}_{\beta}:=\left\{(X,V) :  \left| \Psi^N(\overline{X}^N_t,\overline{V}^N_t)- \overline{\Psi}^N(\overline{X}^N_t, \overline{V}^N_t) \right| _{\infty}> N^{-\beta}\right\},
  \end{eqnarray}
  \begin{eqnarray}\label{ngamma}
 \mathcal{N}_{\gamma}:=\left\{(X,V) : \left|Q^N(\overline{X}^N_t,\overline{V}^N_t)- \overline{Q}^N(\overline{X}^N_t,\overline{V}^N_t)\right|_{\infty}>
 N^{-\gamma}\right\},
 \end{eqnarray}
where $Q^N(\overline{X}^N_t,\overline{V}^N_t)$ and $\overline{Q}^N(\overline{X}^N_t,\overline{V}^N_t)$ are
understood in the sense of
$$ (Q^N(\overline{X}^N_t,\overline{V}^N_t))_i := \frac{1}{N-1}\sum_{i \neq j}q^N(\overline{x}^N_i-\overline{x}^N_j, \overline{v}^N_i-\overline{v}^N_j)$$
and correspondingly
$$ (\overline{Q}^N(\overline{X}^N_t,\overline{V}^N_t))_i := \iint q^N(\overline{x}^N_i-y, \overline{v}^N_i-w)f^N(t,y,w)\,dydw.$$
\end{Definition}

Next, we will see that the measures of both sets $\mathcal{N}_{\beta}$ and $ \mathcal{N}_{\gamma}$
can be arbitrarily small, i.e., the probability of each set tends to 0 as $N$ goes to infinity. We prove the following two lemmas:

 \begin{Lemma}
   There exists  a constant $C < \infty$ such that
   $$\mathbb{P}_0(\mathcal{N}_{\beta}) \le CN^{-(1-4\beta)}.$$
 \end{Lemma}

\textit{Proof.} First, we let the set $\mathcal{N}_{\beta}$ evolve along the
characteristics of the Vlasov equation
$$\mathcal{N}_{\beta,t}:=\left\{(\overline{X}^N_t,\overline{V}^N_t) :  \left|N^{\beta} \Psi^N(\overline{X}^N_t,\overline{V}^N_t)-
N^{\beta}\overline{\Psi}^N(\overline{X}^N_t, \overline{V}^N_t) \right| _{\infty}> 1 \right\}$$
and consider the following fact
$$\mathcal{N}_{\beta,t} \subseteq \bigoplus^{N}_{i=1}\mathcal{N}^i_{\beta,t},$$
where
$$\mathcal{N}^i_{\beta,t}:=\left\{(\overline{x}^N_i,\overline{v}^N_i) :  \left|N^{\beta}\cdot \frac{1}{N-1}\sum_{i \neq j}F^N(\overline{x}^N_i-\overline{x}^N_j,\overline{v}^N_i-\overline{v}^N_j)
-N^{\beta} (F^N*f^N)(t, \overline{x}^N_i,\overline{v}^N_i) \right| _{\infty}> 1 \right\}.$$

\noindent We therefore get
$$\mathbb{P}_t(\mathcal{N}_{\beta,t}) \le \sum_{i=1}^N
\mathbb{P}_t(\mathcal{N}^i_{\beta,t})=N\mathbb{P}_t(\mathcal{N}^1_{\beta,t}),$$
where in the last step we use the symmetry argument in exchanging any two coordinates.

Using Markov inequality gives
\begin{eqnarray}
\label{betaexpectation1}
\mathbb{P}_t(\mathcal{N}^1_{\beta,t}) &\le& \mathbb{E}_t\left[ \left( N^{\beta}\cdot \frac{1}{N-1}\sum^N_{ j =2}F^N(\overline{x}^N_1-\overline{x}^N_j,\overline{v}^N_1-\overline{v}^N_j)
-N^{\beta} (F^N*f^N)(t,\overline{x}^N_1,\overline{v}^N_1)\right)^4 \right]  \nonumber \\
&=& \left(\frac{N^{\beta}}{{N-1}}\right)^4 \mathbb{E}_t\left[ \left( \sum^N_{ j =2}F^N(\overline{x}^N_1-\overline{x}^N_j,\overline{v}^N_1-\overline{v}^N_j)
-(N-1)(F^N*f^N)(t,\overline{x}^N_1,\overline{v}^N_1)\right)^4 \right]. \nonumber \\
\end{eqnarray}

\noindent Let $ \displaystyle h_j:=F^N(\overline{x}^N_1-\overline{x}^N_j,\overline{v}^N_1-\overline{v}^N_j)
-\iint F^N(\overline{x}^N_1-y,\overline{v}^N_1-w)f^N(t,y,w)\,dydw$. Then, each term in
the expectation (\ref{betaexpectation1}) takes the form of $
\prod_{j=2}^Nh_j^{k_j}$ with $\sum_{j=1}^Nk_j=4 $, and more importantly, the expectation assumes the value
of zero whenever there exists a $j$ such that $k_j=1$. This can be easily
verified by integrating over the $j$-th variable first or, in other words, by
acknowledging the fact that $\forall \, j=2, \ldots, N$, there holds
$$ \mathbb{E}_t \left[F^N(\overline{x}^N_1-\overline{x}^N_j,\overline{v}^N_1-\overline{v}^N_j)
-\iint F^N(\overline{x}^N_1-y,\overline{v}^N_1-w)f^N(t,y,w)\,dydw \right] =0.
$$
Then, we can simplify the estimate (\ref{betaexpectation1}) to
\begin{eqnarray*}
\mathbb{P}_t(\mathcal{N}^1_{\beta,t}) \le \left(\frac{N^{\beta}}{{N-1}}\right)^4 \mathbb{E}_t\left[\, \sum_{j=2}^N h_j^4
+\sum^N_{2 \le m<n}\begin{pmatrix} 4 \\ 2 \end{pmatrix} h_m^2h_n^2 \right].
\end{eqnarray*}

\noindent Since $F^N$ is bounded and $||f^N||_1=1$, we thus have  for any fixed $j$
$$
 |h_j| \le |F^N(\overline{x}^N_1-\overline{x}^N_j,\overline{v}^N_1-\overline{v}^N_j)| + \iint |F^N(\overline{x}^N_1-
 y,\overline{v}^N_1-w)|f^N(t,y,w)\,dydw \le C.
$$
Therefore $|h_j|$ is bounded to any power and we obtian
$$  \mathbb{E}_t \left[ h_m^2h_n^2 \right] \le C \quad \hbox{and} \quad \mathbb{E}_t \left[ h_j^4 \right] \le C$$

%

\noindent and consequently
$$\mathbb{P}_t(\mathcal{N}^1_{\beta,t}) \le \left(\frac{N^{\beta}}{{N-1}}\right)^4 \cdot \left(  C \cdot (N-1)+
C\cdot \frac{(N-1)(N-2)}{2}  \right) \le C\cdot N^{-(2-4\beta)}.$$

\noindent By noticing the fact that
$$\mathbb{P}_0(\mathcal{N}_{\beta})=\mathbb{P}_t\left(\mathcal{N}_{\beta,t}\right)\le N\mathbb{P}_t(\mathcal{N}^1_{\beta,t})
\le N \cdot C \cdot N^{-(2-4\beta)}=C\cdot N^{-(1-4\beta)},$$
we obtain the desired result.
\qed\\

  \noindent In fact, this result holds for any $\beta$ if we change accordingly the
  power in the proof to be another even number (depending on $\beta$) greater than four. \\

Due to the singularity of $\nabla_xF$, which is also the motivation for the cut-off, we exploit
a slightly different technique as in Lemma 4.1 to prove

 \begin{Lemma}
   There exists a constant $C < \infty$ such that
   $$\mathbb{P}_0(\mathcal{N}_{\gamma}) \le C \cdot \widetilde{r}(N),$$
where $\widetilde{r}(N)$ is the convergence rate, which is
 $N^{-(1-4\gamma)}\ln^2N$ if $f^N \in L^{\infty}(\R^2 \times \R^2)$
 or $N^{-(1-4\theta-4\gamma)}$ otherwise.
 \end{Lemma}

\textit{Proof.} Let the set $\mathcal{N}_{\gamma}$ evolve along the
characteristics of the Vlasov equation
$$\mathcal{N}_{\gamma,t}:=\left\{(\overline{X}^N_t,\overline{V}^N_t) :  \left|N^{\gamma} Q^N((\overline{X}^N_t,\overline{V}^N_t)-
N^{\gamma}\overline{Q}^N((\overline{X}^N_t,\overline{V}^N_t) \right| _{\infty}> 1 \right\}$$
and consider the fact
$$\mathcal{N}_{\gamma,t} \subseteq \bigoplus^{N}_{i=1}\mathcal{N}^i_{\gamma,t},$$
where
$$\mathcal{N}^i_{\gamma,t}:=\left\{(\overline{x}^N_i,\overline{v}^N_i) :  \left|N^{\gamma}\cdot
\frac{1}{N-1}\sum_{i \neq j}q^N(\overline{x}^N_i-\overline{x}^N_j,\overline{v}^N_i-\overline{v}^N_j)
-N^{\gamma} (q^N*f^N)(t,\overline{x}^N_i,\overline{v}^N_i) \right| _{\infty}> 1 \right\}.$$

\noindent Due to the symmetry in exchanging any two coordinates, we get
$$\mathbb{P}_t(\mathcal{N}_{\gamma,t}) \le \sum_{i=1}^N
\mathbb{P}_t(\mathcal{N}^i_{\gamma,t})=N\mathbb{P}_t(\mathcal{N}^1_{\gamma,t}).$$

\noindent Using Markov inequality gives
\begin{eqnarray}
\label{gammaexpectation1}
\mathbb{P}_t(\mathcal{N}^1_{\gamma,t}) &\le& \mathbb{E}_t\left[ \left( N^{\gamma}\cdot
\frac{1}{N-1}\sum^N_{ j =2}q^N(\overline{x}^N_1-\overline{x}^N_j,\overline{v}^N_1-\overline{v}^N_j)
-N^{\gamma} (q^N*f^N)(t,\overline{x}^N_1,\overline{v}^N_1)\right)^4 \right]  \nonumber \\
&=& \left(\frac{N^{\gamma}}{{N-1}}\right)^4 \mathbb{E}_t\left[ \left( \sum^N_{ j =2}q^N(\overline{x}^N_1-\overline{x}^N_j,\overline{v}^N_1-\overline{v}^N_j)
-(N-1)(q^N*f^N)(t,\overline{x}^N_1,\overline{v}^N_1)\right)^4 \right]. \nonumber \\
\end{eqnarray}

\noindent In order to avoid redundant complexity, we borrow the notation from Lemma 4.1
and also define $\displaystyle h_j:=q^N(\overline{x}^N_1-\overline{x}^N_j,\overline{v}^N_1-\overline{v}^N_j)
-\iint q^N(\overline{x}^N_1-y,\overline{v}^N_1-w)f^N(t,y,w)\,dydw$. With the same
argument as in Lemma 4.1, we simplify the estimate (\ref{gammaexpectation1}) to
$$\mathbb{P}_t(\mathcal{N}^1_{\gamma,t}) \le  \left(\frac{N^{\gamma}}{{N-1}}\right)^4 \mathbb{E}_t\left[\, \sum_{j=2}^N h_j^4
+\sum^N_{2 \le m<n}6h_m^2h_n^2 \right].$$

\noindent On one hand, due to the cut-off, it is clear that
$$
||q^N||_{\infty} \le C \cdot N^{\theta}.
$$

\noindent On the other hand, by taking out the $L^{\infty}$-norm of $q^N$ and using the
integrability of $f^N$, we achieve
$$
\left|\iint q^N(\overline{x}^N_1-y,\overline{v}^N_1-w)f^N(t,y,w)\,dydw\right| \le C \cdot N^{\theta}.
$$

\noindent Therefore $|h_j|$ is bounded by $C \cdot N^{\theta}$ and it is now obvious
to see that
$$ \mathbb{E}_t \left[ h_j^4 \right] \le C\cdot N^{4\theta} \quad \hbox{and} \quad \mathbb{E}_t \left[ h_m^2h_n^2 \right] \le C \cdot N^{4\theta}$$

\noindent and consequently
$$\mathbb{P}_t(\mathcal{N}^1_{\gamma,t}) \le C\cdot \left(\frac{N^{\gamma}}{{N-1}}\right)^4
\left(N^{4\theta}\cdot (N-1)+N^{4\theta}\cdot \frac{ (N-1)(N-2)}{2} \right) \le C\cdot
N^{-(2-4\theta-4\gamma)}.$$

\noindent By noticing the fact that
$$\mathbb{P}_0(\mathcal{N}_{\gamma})=\mathbb{P}_t\left(\mathcal{N}_{\gamma,t}\right)\le N\mathbb{P}_t(\mathcal{N}^1_{\gamma,t})
\le C\cdot N^{-(1-4\theta-4\gamma)},$$
we complete the first part of the lemma. \\

Furthermore , if  $f^N(t,x,v) \in L^{\infty}((0,\infty); L^1(\R^2\times \R^2)) \cap
  L^{\infty}( (0,\infty);L^{\infty}(\R^2\times \R^2))$, by applying the inequality $\mathbb{E}\left[(X-\mathbb{E}[X])^2\right] \le \mathbb{E}[X^2]$ for
  any random variable $X$, we have for any fixed $j$
  \begin{eqnarray*}
&& \mathbb{E}_t \left[\left(q^N(\overline{x}^N_1-\overline{x}^N_j,\overline{v}^N_1-\overline{v}^N_j)
-\iint q^N(\overline{x}^N_1-y,\overline{v}^N_1-w)f^N(t,y,w)\,dydw\right)^2 \right]\\
&\le& \mathbb{E}_t \left[ \left(q^N(\overline{x}^N_1-\overline{x}^N_j,\overline{v}^N_1-\overline{v}^N_j)\right)^2  \right]
\\
&\le& \iint \left( \iint_{|z-y|<N^{-\theta}} \left(C\cdot N^{\theta}+C\right)^2
f^N(t,y, w)\,dydw \right) f^N(t,z, u)\,dz du \\
&& +  \iint \left( \iint_{|z-y|\ge N^{-\theta}} \left(C\cdot \frac{1}{|z-y|}\right)^2
f^N(t,y, w)\,dydw \right) f^N(t,z, u)\,dz du.
\end{eqnarray*}

\noindent We take out the $L^{\infty}$-norm of $f^N$ in both terms. The integral
left inside the first term is bounded by a constant while in the second term the integral
can be estimated by
$$
\iint_{|z-y|\ge N^{-\theta}} \left(C\cdot \frac{1}{|z-y|}\right)^2
f^N(t,y, w)\,dydw \le C +
2\pi \theta \ln N \le C \cdot \ln N,
$$
where we use that $q^N$ has compact support.
Therefore for any fixed $j$
$$
 \mathbb{E}_t \left[ h_j^4 \right] \le ||h_j||_{\infty}^2  \mathbb{E}_t \left[ h_j^2 \right]
 \le C  \cdot N^{2\theta}  \ln N,
$$
$$
 \mathbb{E}_t \left[ h_m^2h_n^2 \right] \le C\cdot \ln^2 N.
 $$
 Consequently
 \begin{eqnarray*}
\mathbb{P}_t(\mathcal{N}^1_{\gamma,t}) &\le& C\cdot \left(\frac{N^{\gamma}}{{N-1}}\right)^4
\left(N^{2\theta}  \ln N \cdot (N-1)+\ln^2 N\cdot \frac{ (N-1)(N-2)}{2} \right) \\
&\le& C\cdot N^{-(2-4\gamma)}\ln^2N.
\end{eqnarray*}

\noindent Thus it holds that
$$\mathbb{P}_0(\mathcal{N}_{\gamma})=\mathbb{P}_t\left(\mathcal{N}_{\gamma,t}\right)\le N\mathbb{P}_t(\mathcal{N}^1_{\gamma,t})
\le C\cdot N^{-(1-4\gamma)}\ln^2N.$$

\qed \\

\begin{Lemma}
  Let $\mathcal{N}_{\al}$, $ \mathcal{N}_{\beta}$, $ \mathcal{N}_{\gamma}$ be defined
  as in (\ref{nalpha})-(\ref{ngamma}). Suppose that $f^N(t,x,v)$ satisfies Assumption 4.1(a) and Assumption 4.1(b) holds for $G(x,v)$. Then there exists a constant $C<\infty$ such that
  $$ \left|\Big(V^N_t, \Psi^N(X^N_t,V^N_t)+\Gamma(X^N_t,V^N_t)\Big)-\Big(\overline{V}^N_t, \overline{\Psi}^N(\overline{X}^N_t, \overline{V}^N_t)+\Gamma(\overline{X}^N_t,\overline{V}^N_t)\Big)
  \right| _{\infty} \le C S_t(X,V)N^{-\al}+N^{-\beta}$$
  for all initial data $(X,V) \in (\mathcal{N}_{\al} \cup \mathcal{N}_{\beta} \cup
  \mathcal{N}_{\gamma})^c$.
\end{Lemma}

\textit{Proof.} Applying triangle inequality gives
\begin{eqnarray*}
&& \left|\Big(V^N_t, \Psi^N(X^N_t,V^N_t)+\Gamma(X^N_t,V^N_t)\Big)-\Big(\overline{V}^N_t, \overline{\Psi}^N(\overline{X}^N_t, \overline{V}^N_t)+\Gamma(\overline{X}^N_t,\overline{V}^N_t)\Big)
  \right| _{\infty}\\
  &\le& \left|V^N_t-\overline{V}^N_t\right| _{\infty}+\left|\Psi^N(X^N_t,V^N_t)-\overline{\Psi}^N(\overline{X}^N_t, \overline{V}^N_t)\right| _{\infty}+
 \left|\Gamma(X^N_t,V^N_t)-\Gamma(\overline{X}^N_t,\overline{V}^N_t)\right| _{\infty} \\
&\le& \left|V^N_t-\overline{V}^N_t\right| _{\infty} +\left|\Psi^N(X^N_t,V^N_t)-\Psi^N(\overline{X}^N_t,\overline{V}^N_t)\right|_{\infty}  \\
&& +\left|\Psi^N(\overline{X}^N_t,\overline{V}^N_t)- \overline{\Psi}^N(\overline{X}^N_t, \overline{V}^N_t)\right|_{\infty}
 + \left|\Gamma(X^N_t,V^N_t)-\Gamma(\overline{X}^N_t,\overline{V}^N_t)\right|
 _{\infty}\\
 &=:& |I_1|+|I_2|+|I_3|+|I_4|.
\end{eqnarray*}
Next, we estimate term by term.
\begin{itemize}
  \item Since $(X,V) \notin \mathcal{N}_{\al}$,
  $$|I_1|:=\left|V^N_t-\overline{V}^N_t\right| _{\infty} \le S_t(X,V)N^{-\al}.$$

\item  With the help of $q^N$ which is defined in Remark 3.1 and the fact that $F^N$ is
Lipschitz continuous in $v$, we obtain
 \begin{eqnarray}
   \label{I2}
   && \left|\frac{1}{N-1}\sum_{i \neq j}F^N(x^N_i-x^N_j,v^N_i-v^N_j)-\frac{1}{N-1}\sum_{i \neq j}F^N(\overline{x}^N_i-\overline{x}^N_j,
  \overline{v}^N_i-\overline{v}^N_j) \right|  \nonumber \\
   &\le& \frac{1}{N-1}\sum_{i \neq j} \big|q^N(\overline{x}^N_i-\overline{x}^N_j,
  \overline{v}^N_i-\overline{v}^N_j)\big| \left( 2|x^N_i-\overline{x}^N_i|+2|v^N_i-\overline{v}^N_i|
   \right).
 \end{eqnarray}
Since $(X,V) \notin \mathcal{N}_{\al}$, it follows in particular for any $1 \le i \le 
N$ that
 $$ |x^N_i-\overline{x}^N_i|\le N^{-\al} \quad \hbox{and} \quad |v^N_i-\overline{v}^N_i|\le N^{-\al}.$$

 So together with (\ref{I2}), we have
$$\left|\big(\Psi^N(X^N_t,V^N_t)\big)_i-\big(\Psi^N(\overline{X}^N_t,\overline{V}^N_t)\big)_i\right| \le 4 \left|\big(Q^N(\overline{X}^N_t,\overline{V}^N_t)\big)_i\right|N^{-\al}. $$

On the other hand, because $(X,V) \notin \mathcal{N}_{\gamma}$, it follows
$$\left|\big(Q^N(\overline{X}^N_t,\overline{V}^N_t)\big)_i\right| \le ||q^N*f^N||_{\infty}+N^{-\gamma} <C$$
and thus
$$|I_2|:=\left|\Psi^N(X^N_t,V^N_t)-\Psi^N(\overline{X}^N_t,\overline{V}^N_t)\right|_{\infty} \le CS_t(X,V)N^{-\al}.$$

  \item Since $(X,V) \notin \mathcal{N}_{\beta}$, it follows directly
 $$|I_3|:=\left|\Psi^N(\overline{X}^N_t,\overline{V}^N_t)- \overline{\Psi}^N(\overline{X}^N_t, \overline{V}^N_t)\right|_{\infty}
 \le N^{-\beta}.$$

 \item Since $G(x,v)$ under Assumption 4.1(b) is Lipschitz continuous, we have for each $1 \le i \le N$
 and $(x^N_i,v^N_i) = \big((X^N_t, V^N_t)\big)_i$, $(\overline{x}^N_i,\overline{v}^N_i) =\big ((\overline{X}^N_t,
 \overline{V}^N_t)\big)_i$
 $$\left|G(x^N_i,v^N_i)-G(\overline{x}^N_i,\overline{v}^N_i)\right| \le L\, \left|(x^N_i,v^N_i)-(\overline{x}^N_i,\overline{v}^N_i)\right|.$$
 Together with the fact that $(X,V) \notin \mathcal{N}_{\al}$, there holds
 $$|I_4|:=\left|\Gamma(X^N_t,V^N_t)-\Gamma(\overline{X}^N_t,\overline{V}^N_t)\right|
 _{\infty} \le LS_t(X,V)N^{-\al}.$$

\end{itemize}

 \noindent Combining all the four terms, we end up with
  $$ \left|\big(V^N_t, \Psi^N(X^N_t,V^N_t)+\Gamma(X^N_t,V^N_t)\big)-\big(\overline{V}^N_t, \overline{\Psi}^N(\overline{X}^N_t, \overline{V}^N_t)+\Gamma(\overline{X}^N_t,\overline{V}^N_t)\big)
  \right| _{\infty} \le C S_t(X,V)N^{-\al}+N^{-\beta}$$
  for all $(X,V) \in (\mathcal{N}_{\al} \cup \mathcal{N}_{\beta} \cup
  \mathcal{N}_{\gamma})^c$.
 \qed\\

\noindent Using Lemmas 4.1 - 4.3 we can now prove Theorem 4.1 and Theorem 4.2:

\subsubsection*{Proof of Theorem 4.1}
From the definition of the Newtonian flow (\ref{NF}) and the characteristics of the Vlasov equation
(\ref{VF}), we know
\begin{eqnarray*}
  (X^N_{t+dt}, V^N_{t+dt}) &=& (X^N_t, V^N_t)+(V^N_t, \Psi^N(X^N_t,
  V^N_t)+\Gamma(X^N_t,V^N_t))dt+o(dt),\\
  (\overline{X}^N_{t+dt}, \overline{V}^N_{t+dt}) &=& (\overline{X}^N_t, \overline{V}^N_t)+(\overline{V}^N_t, \overline{\Psi}^N(\overline{X}^N_t,
  \overline{V}^N_t)+\Gamma(\overline{X}^N_t,\overline{V}^N_t))dt+o(dt).
\end{eqnarray*}
Thus
$$
  \left| (X^N_{t+dt}, V^N_{t+dt}) -  (\overline{X}^N_{t+dt},
  \overline{V}^N_{t+dt})\right|_{\infty} \le \left|(X^N_t, V^N_t)- (\overline{X}^N_t, \overline{V}^N_t)
  \right|_{\infty} $$
$$  +\left|\Big(V^N_t, \Psi^N(X^N_t, V^N_t)+\Gamma(X^N_t,V^N_t)\Big)-\Big(\overline{V}^N_t, \overline{\Psi}^N(\overline{X}^N_t,
  \overline{V}^N_t)+\Gamma(\overline{X}^N_t,\overline{V}^N_t)\Big)\right|_{\infty}dt+o(dt),
$$
i.e.,
$$S_{t+dt}-S_t \le \left|\Big(V^N_t, \Psi^N(X^N_t, V^N_t)+\Gamma(X^N_t,V^N_t)\Big)-\Big(\overline{V}^N_t, \overline{\Psi}^N(\overline{X}^N_t,
  \overline{V}^N_t)+\Gamma(\overline{X}^N_t,\overline{V}^N_t)\Big)\right|_{\infty} N^{\al}dt+o(dt)
  $$

\noindent Taking the expectation over both sides yields
\begin{eqnarray*}
 \mathbb{E}_0\left[\,S_{t+dt}-S_t\right]  &=&  \mathbb{E}_0\left[S_{t+dt}-S_t \,|\, \mathcal{N}_{\al}\,\right]
  +\mathbb{E}_0\left[\,S_{t+dt}-S_t \,|\, \mathcal{N}_{\al}^c\,\right]
  \\
  &\le&  \mathbb{E}_0\left[\,S_{t+dt}-S_t \,|\, (\mathcal{N}_{\beta} \cup \mathcal{N}_{\gamma})
  \setminus \mathcal{N}_{\al}\,\right]
  +\mathbb{E}_0\left[\,S_{t+dt}-S_t \,|\, (\mathcal{N}_{\al} \cup \mathcal{N}_{\beta}
 \cup \mathcal{N}_{\gamma})^c\,\right] \\
  & \le &  \mathbb{E}_0 \left[ \, \left|V^N_t-\overline{V}^N_t \right|_{\infty}\,\Big|\,
  (\mathcal{N}_{\beta} \cup \mathcal{N}_{\gamma}) \setminus \mathcal{N}_{\al}\, \right]N^{\al}dt \\
  && +\, \mathbb{E}_0 \left[ \, \left|\Psi^N(X^N_t, V^N_t)- \overline{\Psi}^N(\overline{X}^N_t,
  \overline{V}^N_t) \right|_{\infty}\,\Big|\, (\mathcal{N}_{\beta} \cup \mathcal{N}_{\gamma}) \setminus \mathcal{N}_{\al}\, \right]N^{\al}dt \\
  && +\, \mathbb{E}_0 \left[ \, \left|\Gamma(X^N_t,V^N_t)-\Gamma(\overline{X}^N_t,\overline{V}^N_t)\right|_{\infty}\,\Big|\,
(\mathcal{N}_{\beta} \cup \mathcal{N}_{\gamma}) \setminus \mathcal{N}_{\al}\, \right]N^{\al}dt\\
&& +\, \mathbb{E}_0\left[\,S_{t+dt}-S_t \,|\, (\mathcal{N}_{\al} \cup \mathcal{N}_{\beta}
 \cup \mathcal{N}_{\gamma})^c\,\right]+o(dt)\\
 &=:& J_1+J_2+J_3+J_4+o(dt),
   \end{eqnarray*}
   where in the second step we use $\mathbb{E}_0(S_{t+dt}-S_{t} \,|\,\mathcal{N}_{\al}) \le 0$  and decompose
   the set $\mathcal{N}^c_{\al} \, $ into
  $(\mathcal{N}_{\beta} \cup \mathcal{N}_{\gamma}) \setminus \mathcal{N}_{\al} \,$ and
  $ (\mathcal{N}_{\al} \cup \mathcal{N}_{\beta} \cup
  \mathcal{N}_{\gamma})^c\,$.\\

\noindent
Since $(X,V) \notin \mathcal{N}_{\al}$, it follows
 \begin{eqnarray*}
J_1 &=& \mathbb{E}_0 \left[ \, \left|V^N_t-\overline{V}^N_t \right|_{\infty}\,\Big|\,
  (\mathcal{N}_{\beta} \cup \mathcal{N}_{\gamma}) \setminus \mathcal{N}_{\al}\, \right]N^{\al}dt
  \\
&\le&
\big(\mathbb{P}_0(\mathcal{N}_{\beta})+\mathbb{P}_0(\mathcal{N}_{\gamma})\big)dt.
  \end{eqnarray*}

\noindent Due to the definition of $\Psi^N$, $\overline{\Psi}^N$, $\Gamma$ as well as the
  boundedness of $F^N$, we obtain

$$ J_2 \le  \left( ||F^N||_{\infty}+ ||F^N*f ||_{\infty} \right) \big(\mathbb{P}_0(\mathcal{N}_{\beta})+
\mathbb{P}_0(\mathcal{N}_{\gamma})\big)N^{\al}dt,$$
$$ J_3  \le 2 ||G||_{\infty}  \big(\mathbb{P}_0(\mathcal{N}_{\beta})
+\mathbb{P}_0(\mathcal{N}_{\gamma})\big)N^{\al}dt.
$$

\noindent Thanks to Lemma 4.1 and Lemma 4.2, we get
\begin{eqnarray*}
J_1+J_2+J_3 &=& \left[ N^{-\al} + C \,   \right] \big(\mathbb{P}_0(\mathcal{N}_{\beta})+\mathbb{P}_0(\mathcal{N}_{\gamma})\big)N^{\al}dt
\\
&\le& C \cdot \max \{ \widetilde{r}(N), N^{-(1-4\beta)} \} N^{\al} dt
 \end{eqnarray*}
 where $\widetilde{r}(N)$ is the convergence rate, which is
 $N^{-(1-4\gamma)}\ln^2N$ if $f^N \in L^{\infty}(\R^2 \times \R^2)$
 or $N^{-(1-4\theta-4\gamma)}$ otherwise.
On the other hand, Lemma 4.3 states that
   \begin{eqnarray*}
     J_4 &=& \mathbb{E}_0\left[\,S_{t+dt}-S_t \,|\, (\mathcal{N}_{\al} \cup \mathcal{N}_{\beta}
 \cup \mathcal{N}_{\gamma})^c\,\right] \\
 &\le& (C \cdot \mathbb{E}_0\left[S_t\right]N^{-\al}+N^{-\beta})\cdot N^{\al}dt +o(dt)\\
 &=& C \cdot \mathbb{E}_0\left[S_t\right]dt+N^{\al-\beta}dt+o(dt).
   \end{eqnarray*}
Therefore, we can determine the estimate
\begin{eqnarray*}
\mathbb{E}_0\left[S_{t+dt}\right]-\mathbb{E}_0[S_{t}] &\le& \mathbb{E}_0\left[S_{t+dt}-S_t \right] \\
&\le& C \cdot \mathbb{E}_0\left[S_t\right]dt
+C \cdot \max\{\widetilde{r}(N) N^{\al},N^{-(1-\al-4\beta}), N^{\al-\beta}\}dt+o(dt).
\end{eqnarray*}
Equivalently, we have
$$ \frac{d}{dt}\,\mathbb{E}_0[S_{t}] \le C\cdot \mathbb{E}_0[S_{t}] +C\cdot \max\{\widetilde{r}(N) N^{\al},N^{-(1-\al-4\beta}), N^{\al-\beta}\}.$$

\noindent Gronwall's inequality yields
$$ \mathbb{E}_0\left[S_t\right] \le e^{Ct}\cdot \max\{\widetilde{r}(N) N^{\al},N^{-(1-\al-4\beta}), N^{\al-\beta}\}.$$
The proof is completed by the following Markov inequality
$$\mathbb{P}_0 \left(\sup_{0\le s \le t}\Big|(X^N_s,V^N_s)-(\overline{X}^N_s, \overline{V}^N_s)\Big|_{\infty}> N^{-\al} \right) =
\mathbb{P}_0 (S_t=1) \le  \mathbb{E}_0\left[S_t\right]. $$
\qed

 \subsubsection*{Proof of Theorem 4.2}
Let $N \in \N$ and
  $$W_t:=\sup_{(X,V)\in \R^{4N}} \left|(\overline{X}^N_t, \overline{V}^N_t)-(\overline{X}_t, \overline{V}_t) \right|.$$
With the same argument as in the proof of Theorem 4.1, it is not difficult to deduce
 $$W_{t+dt}-W_t \le \underbrace{ \left|\big(\overline{V}^N_t, \overline{\Psi}^N(\overline{X}^N_t, \overline{V}^N_t)
    +\Gamma(\overline{X}^N_t,\overline{V}^N_t)\big) -\big(\overline{V}_t, \overline{\Psi}(\overline{X}_t,\overline{V}_t)
    +\Gamma(\overline{X}_t,\overline{V}_t)\big)\right|_{\infty}}_{=:D}dt+o(dt).
  $$

Furthermore, with the Lipschitz continuity of $G(x,v)$, we get
  \begin{eqnarray*}
    D  &\le&  W_t +  \left|\overline{\Psi}^N(\overline{X}^N_t, \overline{V}^N_t)-\overline{\Psi}(\overline{X}_t,\overline{V}_t)\right|_{\infty}+
   \left|\Gamma(\overline{X}^N_t,\overline{V}^N_t)-\Gamma(\overline{X}_t,\overline{V}_t)\right|_{\infty} \\
     & \le & C \cdot W_t +  \left|\overline{\Psi}^N(\overline{X}^N_t,
     \overline{V}^N_t)-\overline{\Psi}(\overline{X}_t,\overline{V}_t)\right|_{\infty}\\
      & \le & C \cdot W_t +  \sup_{1 \le i \le N}\Big|F^N*f^N(\overline{x}^N_i,\overline{v}^N_i)-F*f(\overline{x}_i,\overline{v}_i)\Big|\\
      &\le &  C \cdot W_t +  \sup_{1 \le i \le N}\Big|F^N*f^N(\overline{x}^N_i,\overline{v}^N_i)-F^N*f^N(\overline{x}_i,\overline{v}_i)\Big|\\
      && + \sup_{1 \le i \le N}\Big|F^N*f^N(\overline{x}_i,\overline{v}_i)-F^N*f(\overline{x}_i,\overline{v}_i)\Big| \\
      &&
      +\sup_{1 \le i \le N}\Big|F^N*f(\overline{x}_i,\overline{v}_i)-F*f(\overline{x}_i,\overline{v}_i)\Big|.
     \end{eqnarray*}


     \noindent By using the integrability of $\nabla F^N$ , we estimate the second term
      by
          $$ \sup_{1 \le i \le N}\Big|F^N*f^N(\overline{x}^N_i,\overline{v}^N_i)-F^N*f^N(\overline{x}_i,\overline{v}_i)\Big|
      \le ||\nabla F^N ||_1||f^N||_{\infty}W_t \le C\cdot W_t.$$


      Due to the integrability of $\nabla f_0$, the third term can be  controlled by
            $$\sup_{1 \le i \le N}\Big|F^N*f^N(\overline{x}_i,\overline{v}_i)-F^N*f(\overline{x}_i,\overline{v}_i)\Big| \le ||F^N||_{\infty}||f^N-f||_{1}
      \le C ||\nabla f_0||_1 W_t.$$

			In the estimates above, the reversibility of both particle trajectories is used.
      The last term is straightforward to estimate
            $$\sup_{1 \le i \le N}\Big|F^N*f(\overline{x}_i,\overline{v}_i)-F*f(\overline{x}_i,\overline{v}_i)\Big|
      \le ||f||_{\infty}||F^N-F||_1 \le C \cdot N^{-\theta}.$$

      Therefore we arrive at
$$
       W_{t+dt}-W_t  \le  \big(C \cdot W_t +  C \cdot N^{-\theta} \big)dt+o(dt),
$$
or equivalently
$$\frac{d}{dt}W_t \le C \cdot W_t +  C \cdot N^{-\theta}.$$
 Gronwall's inequality gives
 $$W_t \le  C \cdot N^{-\theta}.$$
 Together with Theorem 4.1, we complete the proof.
 \qed

\section{Propagation of Chaos} \label{sec:prog}

We can clearly see as the direct byproduct of the results stated above that chaos indeed propagates, which means the convergence of
the one particle marginals of the $N$-particle system to the solution of the
Vlasov equation in the sense of bounded Lipschitz distance. We illustrate the
propagation of chaos also in two steps by using the Vlasov flow with cut-off as
an intermediate tool. We present the result in full detail under the conditions of Theorem 4.1.

\begin{Definition}
  For any two probability densities $\mu$, $\nu : \R^{4} \to \R^+$, the bounded
  Lipschitz distance is defined by
  $$d_L(\mu,\nu):= \sup_{g \in \mathcal{L}}\left| \int \big( \mu (x,v)-\nu (x,v) \big) g(x,v) \,dxdv \right|,$$
  where
  $\mathcal{L}:=\{g  :  ||g||_{\infty}=||g||_L=1 \}$ and $||g||_L$ denotes the
  global Lipschitz constant of $g$.
\end{Definition}

In order to simplify the notation, we also introduce hereafter $(x_{i,-t},v_{i,-t})$ and $( \overline{x}_{i,-t},\overline{v}_{i,-t})$
to be the position and velocity of the $i$-th particle at initial time, which
evolves according to the Newtonian and Vlasov flow with cut-off starting from $(x_i,v_i)$ at time $t$, respectively.

\begin{Theorem}
  Let $f^N_t : \R \times \R^4 \to \R^+$ be the solution to \eqref{vlasovN}, $\mu_t : \R \times \R^{4N} \to \R^+$ be the
  $N$-particle density of the Newtonian flow and the one-particle marginals $\mu_t^{(1)}$ be given by
  $$ \mu_t^{(1)}(x_1,v_1):=\int \mu_t  (x_1,v_1, \cdots, x_N, v_N) \,dx_2 dv_2 \ldots dx_N dv_N, $$
  where
     $$\mu_t (x_1,v_1, \cdots, x_N, v_N) :=\mu_0 (x_{1,-t},v_{1,-t}, \cdots, x_{N,-t}, v_{N,-t}). $$
  Assume that initially the one particle marginals converges to the initial probability density $f^N_0$ in the sense of bounded Lipschitz distance, i.e.,
  $$ \lim_{N \to \infty} d_L(\mu_0^{(1)}, f_0^N)=0.$$
Then under the conditions of Theorem 4.1, there holds
  $$ \lim_{N \to \infty} d_L(\mu_t^{(1)}, f_t^N)=0.$$
\end{Theorem}

\textit{Proof.} By definition, we have
\begin{eqnarray}
  d_L(\mu_t^{(1)}, f_t^N) &=&  \sup_{g \in \mathcal{L}}\left| \int \big( \mu_t^{(1)}(x_1,v_1)-f_t^N(x_1,v_1)  \big) g(x_1,v_1) \,dx_1 dv_1 \right| \nonumber\\
 &=& \sup_{g \in \mathcal{L}}\Big| \int \big( \mu_t  (x_1,v_1, \cdots, x_N, v_N)- \nonumber \\
 \label{distance}
 && \qquad  \prod_{i=1}^N f_t^N (x_i,v_i) \big) g(x_1,v_1) \,dx_1 dv_1dx_2 dv_2 \ldots dx_N dv_N \Big| .
 \end{eqnarray}

\noindent Since both the Newtonian and Vlasov flow leave the measure invariant, then

\begin{eqnarray*}
 \eqref{distance} &=&  \sup_{g \in \mathcal{L}}\Big| \int \mu_0  (x_1,v_1, \cdots, x_N, v_N)g(x_{1,-t},v_{1,-t})\,dx_1 dv_1 \ldots dx_N dv_N \\
 && \qquad -\int \prod_{i=1}^N f_0^N (x_i,v_i) g(\overline{x}_{1,-t},\overline{v}_{1,-t}) \,dx_1 dv_1 \ldots dx_N dv_N \Big| \nonumber \\
  &\le&  \sup_{g \in \mathcal{L}}\left| \int \mu_0 (x_1, v_1, \cdots, x_N, v_N)\big( g(x_{1,-t},v_{1,-t})- g(\overline{x}_{1,-t},\overline{v}_{1,-t}) \big)\,dx_1 dv_1 \ldots dx_N dv_N \right| \nonumber \\
   &&  +  \sup_{g \in \mathcal{L}}\left| \int \big(\mu_0 (x_1, v_1, \cdots, x_N, v_N)-\prod_{i=1}^N f_0^N (x_i,v_i)\big) g(\overline{x}_{1,-t},\overline{v}_{1,-t}) \,dx_1 dv_1 \ldots dx_N dv_N   \right|  \\
&=:&  M_1+M_2.
\end{eqnarray*}
Further we decompose $M_1$ into $M_{11}+M_{12}$, where
$$ M_{11}=M_1\Big|_{\left\{  \sup_{0\le s \le t}\left|(X^N_s,V^N_s)-(\overline{X}^N_s,
\overline{V}^N_s)\right|_{\infty}> N^{-\al}\right\}}$$
and
$$M_{12}=M_1\Big|_{\left\{ \sup_{0\le s \le t}\left|(X^N_s,V^N_s)-(\overline{X}^N_s,
\overline{V}^N_s)\right|_{\infty}\le N^{-\al}\right\}}.$$
 Under Theorem 4.1, we know
 $$\lim_{N \to \infty} \mathbb{P}_0 \left(\sup_{0\le s \le t}\left|(X^N_s,V^N_s)-(\overline{X}^N_s, \overline{V}^N_s)\right|_{\infty}> N^{-\al} \right)
=0.$$
By using the fact that $||g||_{\infty}=1$ , we thus obtain
\begin{eqnarray*}
  M_{11} &<& 2  \int \mu_0 (x_1, v_1, \cdots, x_N, v_N)\,dx_1 dv_1 \ldots dx_N dv_N \\
  &<& 2 \mathbb{P}_0 \left(\sup_{0\le s \le t}\left|(X^N_s,V^N_s)-(\overline{X}^N_s, \overline{V}^N_s)\right|_{\infty}> N^{-\al} \right) \\
&\to& 0, \quad \hbox{as}\,N \to \infty.
\end{eqnarray*}

\noindent On the other hand, due to the reversibility of both particle trajectories and $||g||_{L}=1$, we have
\begin{eqnarray*}
  M_{12} &<&   \int \mu_0 (x_1, v_1, \cdots, x_N, v_N)\big| (x_{1,-t},v_{1,-t})- (\overline{x}_{1,-t},\overline{v}_{1,-t}) \big|\,dx_1 dv_1 \ldots dx_N dv_N \\
  &=& \mathbb{E}_0\big[ \big| (X_{1,-t},V_{1,-t})- (\overline{X}_{1,-t},\overline{V}_{1,-t}) \big|  \big] \\
  &<& \mathbb{E}_0\left[\sup_{0\le s \le t}\left|(X^N_s,V^N_s)-(\overline{X}^N_s, \overline{V}^N_s)\right|_{\infty}
  \right]\\
  &\to& 0, \quad \hbox{as}\,N \to \infty.
\end{eqnarray*}
In summary, $M_1$ converges to zero as $N$ goes to infinity. Meanwhile it is also clear that $M_2$ tends to zero as $N \to \infty$ due to the assumption on the
initial probability density. Combining all the terms completes the proof.
\qed

\begin{Theorem}
  Let $f_t : \R \times \R^4 \to \R^+$ be the solution to \eqref{VE}, $\mu_t : \R \times \R^{4N} \to \R^+$ be the
  $N$-particle density of the Newtonian flow and the one-particle marginals $\mu_t^{(1)}$ be given by
  $$ \mu_t^{(1)}(x_1,v_1):=\int \mu_t  (x_1,v_1, \cdots, x_N, v_N) \,dx_2 dv_2 \ldots dx_N dv_N, $$
  where
     $$\mu_t (x_1,v_1, \cdots, x_N, v_N) :=\mu_0 (x_{1,-t},v_{1,-t}, \cdots, x_{N,-t}, v_{N,-t}). $$
  Assume that initially the one particle marginals converges to the initial probability density $f_0$ in the sense of bounded Lipschitz distance, i.e.,
  $$ \lim_{N \to \infty} d_L(\mu_0^{(1)}, f_0)=0.$$
Then under the conditions of Theorem 4.2, there holds
  $$ \lim_{N \to \infty} d_L(\mu_t^{(1)}, f_t)=0.$$
\end{Theorem}

\textit{Proof.} By replacing $f_t^N$ with $f_t$ in the proof of Theorem 5.1 and
using the conditions of Theorem 4.2, one will directly get the desired result. But
we emphasize that Theorem 5.2 actually implies the convergence of
the solution of \eqref{vlasovN} to the solution of \eqref{VE} in the sense of
bounded Lipschitz distance.
\qed \\

Note that if the initial one particle marginals converges in a
certain rate to the initial probability density in both theorems above, we can also achieve the
convergence rate for any fixed time $t$.


\section{Summary}
This paper deals with one core problem: how to derive rigorously the kinetic
description of one particle density evolution from the $N$-particle system for $N$
being large enough. Our main results, Theorem 4.1 and Theorem 4.2, state that the
trajectories of both the Newtonian system with cut-off and the characteristics of the Vlasov equation are
close to each other in the probability sense.
Propagation of chaos, as the direct implication of the two theorems, is given in Theorem 5.2.
The existence and uniqueness of the $L^{\infty}( (0,\infty);L^{\infty}(\R^2\times \R^2))$-solution of the Vlasov equation
is left for a future independent work.


\section*{Acknowledgments}
This work was financially supported by the DAAD project “DAAD-PPP VR China”
(Project-ID: 57215936).

%

%

\end{document}